\documentclass[12pt, reqno]{amsart}

\usepackage{amsthm}
\usepackage{amssymb}
\usepackage{amsmath}

\newtheorem{theorem}{Theorem}[section]
\newtheorem*{theorem*}{Theorem}
\newtheorem{lemma}{Lemma}[section]
\newtheorem{proposition}{Proposition}[section]
\newtheorem*{proposition*}{Proposition}

\theoremstyle{definition}
\newtheorem{definition}{Definition}[section]

\newtheorem{remark}{Remark}[section]

\numberwithin{equation}{section}

\newcommand{\1}{{1\hspace{-0.35em}1}}

\newcommand{\R}{{\mathbb R}}
\newcommand{\C}{{\mathbb C}}
\newcommand{\Z}{{\mathbb Z}}

\newcommand{\CXbX}{{\C[X, {\bar X}]}}
\newcommand{\CG}{{\C(G)}}
\newcommand{\CP}{{\C\left({\mathbb P}^{n-1|m}\right)}}
\newcommand{\CS}{{\C\left({\mathbb S}^{2n-1|2m}\right)}}
\newcommand{\CKG}{{\C(K\backslash G)}}
\newcommand{\CKGK}{{\C(K\backslash G/K)}}
\newcommand{\Hom}{{\rm{Hom}}}

\newcommand{\g}{{\mathfrak g}}
\newcommand{\fg}{{\mathfrak g}}
\newcommand{\fk}{{\mathfrak k}}
\newcommand{\fl}{{\mathfrak l}}

\newcommand{\fp}{{\mathfrak p}}
\newcommand{\fq}{{\mathfrak q}}
\newcommand{\fu}{{\mathfrak u}}
\newcommand{\fS}{{\mathfrak S}}
\newcommand{\fgu}{{{\mathfrak g}^{\sigma,\sqrt{i}}}}

\newcommand{\U}{{\rm U}({\mathfrak g})}
\newcommand{\Ul}{{\rm U}({\mathfrak l})}
\newcommand{\URk}{{{\rm U}^{\R}({\mathfrak k})}}
\newcommand{\id}{{\rm id}}

\newcommand{\cJ}{{\mathcal J}}

\newcommand{\bea}{\begin{eqnarray}}
\newcommand{\eea}{\end{eqnarray}}
\newcommand{\be}{\begin{eqnarray*}}
\newcommand{\ee}{\end{eqnarray*}}

\begin{document}

\title[Spherical Functions On homogeneous superspaces]
{Spherical Functions On homogeneous superspaces}
\author[R.B. Zhang]{R.B. Zhang}
\address{School of Mathematics and Statistics,
University of Sydney, Sydney, NSW 2006, Australia}
\email{rzhang@maths.usyd.edu.au}
\author[Y.M. Zou]{Y.M. Zou}
\address{Department of Mathematical Sciences, University
of Wisconsin, Milwaukee, WI 53201, USA} \email{ymzou@uwm.edu}

\begin{abstract}
Homogeneous superspaces arising from the general linear supergroup
are studied within a Hopf algebraic framework. Spherical functions
on homogeneous superspaces are introduced, and the structures of
the superalgebras of the spherical functions on classes of
homogeneous superspaces are described explicitly.
\end{abstract}
\maketitle
%\tableofcontents

\section{Introduction}
We study spherical functions on homogeneous superspaces arising
from the complex general linear supergroup. This is the first part
of our endeavour to develop a theory of spherical functions on Lie
supergroups \cite{Ko, Ma88} and quantum supergroups \cite{Ma89,
Zh98}. The theory of spherical functions on ordinary Lie groups
has long reached its maturity (see, e.g., \cite{Ta}). There also
exists extensive literature on spherical functions \cite{Koo, NM,
MNU, DN, No} on quantum symmetric spaces \cite{KD, No, DN, GZ,
Le}. However, little seems to be known about spherical functions
on Lie supergroups, let alone those on quantum supergroups. On the
other hand, supersymmetry and its quantum analogue have become an
integral part of modern mathematical physics, and have also
permeated many areas of pure mathematics. A good understanding of
spherical functions on Lie supergroups and quantum supergroups
should facilitate practical means for studying the dynamics of
physical systems with classical or quantum supersymmetries.

We choose to work within a Hopf superalgebraic framework to study
homogeneous superspaces, as it can incorporate both the Lie
supergroup (as defined by Kostant \cite{Ko}) and quantum
supergroup \cite{Ma89, Zh98} cases.  Our methodology is similar to
that adopted in the literature on quantum homogeneous spaces
\cite{KD, No, DN, GZ, Le}. The starting point is the universal
enveloping algebra $\U$ of the general linear superalgebra $\fg=
{\mathfrak{gl}}(m|n, \C)$, which is a co-commutative Hopf
superalgebra \cite{MM}. A $\Z_2$-graded subalgebra $\CG$ (see
Definition \ref{function}) of the dual of the universal enveloping
algebra acquires a Hopf superalgebra structure, from which the
general linear supergroup can be re-constructed \cite{SZ02} in a
manner similar to the Tanakan-Krein theory for compact Lie groups.
The universal enveloping algebra admits many Hopf
$\ast$-superalgebra structures, each corresponding to a real form
$\fgu$ (see Section \ref{real} for definition) of $\fg$. Each Hopf
$\ast$-superalgebra structure $\theta$ of $\U$ induces a Hopf
$\ast$-superalgebraic structure on $\CG$. We fix the $\theta$
corresponding to one of the compact real forms of $\fg$ (see
equation (\ref{compact-gl})). Let $\fp\subset\fg$ be a parabolic
subalgebra with Levi factor $\fl$, and let $\fk=\fl\cap\fgu$ be
the real form of $\fl$. Then the $\ast$-subalgebra $\CKG$ of $\CG$
invariant with respect to $\fk$ under the left translation defines
a homogeneous superspace \cite{Ma88} in the spirit of
non-commutative geometry \cite{Co}. We shall call this
superalgebra the superalgebra of functions on the homogeneous
superspace. Next we consider the subspace $\CKGK$ of $\CKG$
consisting of elements that are invariant with respect to $\fk$
under the right translation. It can be shown that $\CKGK$ forms a
$\ast$-superalgebra, which will be referred to as the superalgebra
of spherical functions on the homogeneous superspace.

Our aim in the present paper is to understand the structures of
the superalgebras $\CKG$ and $\CKGK$.  The main results obtained
are Theorem \ref{construction}, Lemma \ref{gCKG} and Lemma
\ref{gCKGK}, which give explicit descriptions of the superalgebra
of functions on the homogeneous superspace and the superalgebra of
spherical functions. In the case of a homogeneous superspace
associated to a maximal rank reductive subgroup of a compact real
form of the the general linear supergroup, the superalgebra of
spherical functions is either the polynomial algebra in one
variable or a quotient thereof (Theorems \ref{projective} and
\ref{maxrank}).

Recall that the space of functions on an ordinary Lie group has
another natural algebraic structure with the multiplication
defined by convolution. In this context, the counter parts of
$\CKG$ and $\CKGK$ form subalgebras under convolution, where the
analogue of $\CKGK$ is the celebrated Hecke algebra \cite{Ta}. The
Hecke algebras associated with Riemannian symmetric spaces are
commutative, and their elements provide the invariant integral
operators acting on functions on the symmetric spaces. It is an
important problem to develop a theory for such Hecke algebras in
the Lie supergroup context, and to investigate properties of
supersymmetric spaces from the viewpoint of Hecke algebras. We
plan to do this in a future publication, as the problem requires
in depth investigations into the analytical theory of Lie
supergroups.

The organization of the paper is as follows. In Section
\ref{Preliminaries} we provide some preliminary material on the
complex general linear superalgebra and its invariant theory. In
Section \ref{functions} we discuss the Hopf superalgebra of
functions on the general linear supergroup, and explain how the
general linear supergroup itself can be extracted from this Hopf
superalgebra \cite{SZ02}. The material in this section is not all
new, but it forms the basis for the study of homogenous
superspaces and spherical functions in later sections. Sections
\ref{sect-space} and \ref{max-rank} contain the main results of
the paper. In Subsection \ref{real} we discuss real forms of the
complex general linear superalgebra and general linear supergroup
from a Hopf algebraic point of view. The material presented here
is largely new, and we believe it to be interesting in its own
right. In Subsection \ref{subsect-space} we explain the notion of
homogeneous superspaces in a Hopf algebraic setting, and in In
Subsection \ref{main} we investigate the superalgebras of
spherical functions on the homogeneous superspaces. In Section
\ref{max-rank} we analyze in detail the superalgebras of spherical
functions on the projective superspace and other symmetric
superspaces arising from maximal rank subgroups of real forms of
the general linear supergroup.

\section{Preliminaries on ${\mathfrak{gl}}(m|n,
\C)$}\label{Preliminaries}
We present some background material on the universal enveloping
superalgebra of the general linear Lie superalgebra, which will be
used later. General references are \cite{Kac77, Sc79}.

We shall work on the complex number field $\C$ for simplicity.
Let $W$ be a superspace, i.e., a $\Z_2$-graded vector space
$W=W_{\bar 0}\oplus W_{\bar 1}$, where $W_{\bar 0}$ and $W_{\bar 1}$ are
the even and odd subspaces, respectively. The elements of
$W_{\bar 0}\cup W_{\bar 1}$ will be called homogeneous. Define a map
$[\ ]: W_{\bar 0}\cup W_{\bar 1} \rightarrow \Z_2$ by $[w] = \alpha$ if
$w \in W_{\alpha}$. (Quite generally, whenever a symbol like $[w]$ appears
in the sequel, it is tacitly assumed that the element $w$ is homogeneous.)
The dual superspace ($\Z_2$-graded dual vector space) of $W$ will be
denoted by $W^*$, and the dual space pairing $W^* \otimes W\rightarrow \C$
by $\langle \; , \: \rangle$.

Denote by $\g$ the Lie superalgebra ${\mathfrak{gl}}(m|n,\C)$. A standard
basis for $\g$ is $\{E_{a b}\,|\, a, b\in{\bf I}\},$ where
${\bf I}=\{1, 2, \ldots , m+n\}.$ The element $E_{a b}$ belongs to
$\g_{\bar 1}$ if $a\le m<b$, or $b\le m<a$, and belongs to $\g_{\bar 0}$
otherwise. For convenience, we define the map
\be
[\ ]: {\bf I}\rightarrow \Z_2 \quad {\text by} \quad
[a]=\left\{\begin{array}{l l}
               \bar{0},  & \mbox{if} \ a\le m, \\
               \bar{1},  & \mbox{if} \ a>m.
              \end{array}\right.
\ee
Then $[E_{a b}] = [a]+[b]$. The supercommutation relations of the Lie
superalgebra are given for the basis elements by
\be
[E_{a b}, \ E_{c d}] &=& E_{a d}\delta_{b c}
                         - (-1)^{([a]-[b])([c]-[d])} E_{c b}\delta_{a d}.
\ee
As usual, we choose the Cartan subalgebra
${\mathfrak h}=\bigoplus_a \C E_{a a}$. Let
$\{\epsilon_a \,|\, a\in{\bf I}\}$ be the basis of ${\mathfrak h}^*$ such
that $\epsilon_a(E_{b b})=\delta_{a b}$. The space ${\mathfrak h}^*$ is
equipped with a bilinear form
$(\; , \: ): {\mathfrak h}^*\times {\mathfrak h}^* \rightarrow \C$ such that
$(\epsilon_a, \epsilon_b)=(-1)^{[a]} \delta_{a b}$. The roots of $\g$ are
$\epsilon_a-\epsilon_b, \:\, a\ne b$, where $\epsilon_a-\epsilon_b$ is even
if $[a]+[b] = \bar{0}$ and odd otherwise. We choose as positive roots the
elements of $\{ \epsilon_a-\epsilon_b \,|\, a< b\}$, and as simple roots
the elements of $\{\epsilon_a-\epsilon_{a+1} \,|\, a<m+n\}$.

The enveloping algebra ${\rm U}({{\mathfrak g}{\mathfrak l}}(m|n, \C))$ of
${{\mathfrak g}{\mathfrak l}}(m|n, \C)$ will be denoted by $\U$.
We shall always regard $\fg$ as embedded in $\U$ in the natural way.
As is well known, $\U$ forms a $\Z_2$-graded cocommutative Hopf algebra
(i.e., a Hopf superalgebra) in the sense of \cite{MM}, with
\begin{tabbing}
comultiplication:\quad  \= $\U\rightarrow \U\otimes \U$,\quad  \=
$\Delta(X)=X\otimes 1 + 1\otimes X$, \quad \= $X\in\g$,\\
counit: \>$\epsilon: \U\rightarrow \C$,  \> $\epsilon(X)=0$, \> $X\in\g$, \\
antipode: \> $S: \U\rightarrow\U$,  \>$S(X)=-X$, \>$X\in\g$. \\
\end{tabbing}
In particular, this Hopf superalgebra structure allows
us to introduce a natural left $\U$-module structure on the dual superspace
$W^*$ of any left $\U$-module $W$, with the $\U$-action given by
\be
\U\otimes W^*&\rightarrow& W^*, \quad
x\otimes {\bar w}\mapsto x{\bar w}, \\
\langle x{\bar w}, \, v\rangle&:=& (-1)^{[x][{\bar w}]} \langle {\bar w},
\, S(x) v\rangle, \quad \forall\,  v\in W. \ee
As it stands, the last equation only makes sense for homogeneous
${\bar w}\in W^*$ and homogeneous $x\in\U$, but it can be extended to all
elements of $W^*$ and $\U$ linearly.

We shall denote by $L_{\lambda}$ the irreducible left $\U$-module with
highest weight $\lambda\in{\mathfrak h}^*$. The module $L_{\lambda}$ is
finite-dimensional if and only if $\lambda$ is dominant
\cite{Kac78, Sc79}, i.e.,
\bea {2(\lambda, \, \epsilon_a-\epsilon_{a+1})}/
{(\epsilon_a-\epsilon_{a+1},\, \epsilon_a-\epsilon_{a+1})} &\in& \Z_+
   \quad \forall\, a\ne m. \label{dominant}\eea
A basic problem in the
representation theory of Lie superalgebras is to understand the weight
space decompositions of the finite dimensional irreducible
representations. However, the problem turned out to be unexpectedly
difficult, resisting solution for some twenty years.
Only a few years ago, Serganova \cite{Se} succeeded in developing
an algorithm to compute formal characters of irreducible representations.

Of particular importance to us here is the contravariant vector module
$V=L_{\epsilon_1}$ of $\g$. It has the standard basis
$\{v_a | a\in{\bf I}\}$ such that $E_{a b} v_c = \delta_{b c} v_a$,
where $v_a$ is even if $a\le m$, and odd otherwise.  The
dual module $V^*$ of $V$ is the covariant vector module with
highest weight $-\epsilon_{m+n}$. It has a basis $\{{\bar v}_a \,|\, a\in{\bf I}\}$
dual to the standard basis of $V$, i.e.,
$\langle {\bar v}_a, \, v_b \rangle= \delta_{a b}. $
The action of $\g$ on $V^*$ is given by
\bea
E_{ab}{\bar v}_c &=& -(-1)^{[a]+[a][b]}\delta_{ac}{\bar v}_b.
\eea
As the antipode of $\U$ is of order two, there is a $\U$-module isomorphism
between $V$ and its double dual $V^{**}:=(V^*)^*$:
\be V \cong V^{**}, && v\mapsto v^{**}, \\
\langle v^{**},\ \bar{w}\rangle &=&(-1)^{[v]}
\langle \bar{w}, \ v\rangle,\quad \forall {\bar w}\in V^*.\ee

\begin{remark}\cite{Zh89}
For all $d>0$, $V^{\otimes d}$ is a semi-simple $\U$-module,
which does not contain any $1$-dimensional submodule.
\end{remark}

Let $\fS_d$ be the symmetric group on $d$ letters. There exists a
natural action $\rho_d$ of $\fS_d$ on $V^{\otimes d}$ defined in
the following way. Let $s_i$ denote the permutation $(i,\ i+1)$. Then
\be &&\rho_d(s_i)\left(v_{a_1}\otimes ... \otimes v_{a_{i-1}}\otimes
v_{a_i}\otimes v_{a_{i+1}} \otimes v_{a_{i+2}}
... \otimes v_{a_d}\right)\\
&&= (-1)^{[a_i][a_{i+1}]}v_{a_1}\otimes ... \otimes
v_{a_{i-1}}\otimes v_{a_{i+1}}\otimes v_{a_i} \otimes v_{a_{i+2}}
... \otimes v_{a_d}. \ee Let us denote by $t^d$ the representation
of $\U$ in $V^{\otimes d}$, and denote by $\C\fS_d$ the group
algebra of $\fS_d$. The following result was first proven by Sergeev
\cite{Sergeev, Sergeev01} (see \cite{BR} for a detailed treatment of
the result).
\begin{theorem} \label{FFT}
The superalgebras $t^d(\U)$ and $\rho_d(\C\fS_d)$ are mutual centralizers in
${\rm End}_\C(V^{\otimes d})$.
\end{theorem}

Let $W$ be a finite dimensional $\U$-module.  Let $\pi: \ \U\rightarrow  \rm{End}_{\C}(W)$
be the $\U$-representation furnished by $W$. Then $\rm{End}_{\C}(W)$ acquires a natural
$\U$-module structure under the action
\be \U\otimes \rm{End}_{\C}(W)&\rightarrow& \rm{End}_{\C}(W), \quad
x\otimes \phi \mapsto {\rm Ad}_x(\phi), \\
{\rm Ad}_x(\phi)&:=& \sum_{(x)} (-1)^{[x_{(2)}][\phi]}
\pi(x_{(1)}) \phi \pi(S(x_{(2)}), \ee
where we have used Sweedler's notation $\Delta(x)=\sum_{(x)}
x_{(1)}\otimes x_{(2)}$ for the co-multiplication of $x\in\U$.
There exists the natural isomorphism $j: W\otimes W^* \cong
\rm{End}_{\C}(W)$ of $\U$-modules defined,  for any
$u\otimes {\bar v}\in W\otimes W^*$ and $w\in W$,  by
\be j(u\otimes {\bar v})(w)&=& \langle {\bar v}, w\rangle u. \ee
For any $\U$-module $M$, we use the notation $(M)^{\U}$ to denote the invariant
submodule $$(M)^{\U}:=\{w\in M | x w =\epsilon(x) w, \ \forall
x\in\U\}.$$ We have
\bea  \left(W\otimes W^*\right)^{\U} &\cong&
\rm{End}_{\U}(W)\nonumber\\
&:= &\{ \phi\in {\rm End}_{\C}(W) \mid {\rm
Ad}_x(\phi)=\epsilon(x)\phi, \forall x\in\U\}. \eea Consider
$V^{\otimes k}\otimes (V^*)^{\otimes\ell}$ as a $\U$-module, where
the $\U$-action is defined by using the co-multiplication. The
element $Z=\sum_{a} E_{a a}$ acts on $V^{\otimes k}\otimes
(V^*)^{\otimes\ell}$ by $(k-\ell)\id$. This immediately shows that
\bea \left(V^{\otimes k}\otimes (V^*)^{\otimes\ell}\right)^{\U}
&=&\{0\} , \quad \mbox{if}\ k\ne \ell. \eea As $(V^{\otimes
d})^*\cong(V^*)^{\otimes d}$, we have the $\U$-module isomorphism
$$j: \ V^{\otimes d}\otimes (V^*)^{\otimes d} \rightarrow {\rm
End}_\C(V^{\otimes d}).$$ It follows from Theorem \ref{FFT} that
the even subspace of $\left(V^{\otimes d}\otimes (V^*)^{\otimes
d}\right)^{\U}$ is isomorphic to $j^{-1}\circ\rho_d(\C\fS_d).$ Let
$\fg_{\bar 0}={\mathfrak{gl}}(m)\oplus {\mathfrak{gl}}(n)$ be the
maximal even subalgebra of $\fg$. Both $V$ and $V^*$ naturally
restrict to $\fg_{\bar 0}$-modules. By using Weyl's first
fundamental theorem for the invariant theory of the general linear
group \cite{Ho}, we easily prove that $\left(V^{\otimes d}\otimes
(V^*)^{\otimes d}\right)^{\rm{U}(\fg_{\bar 0})}$ is contained in
the even subspace of $V^{\otimes d}\otimes (V^*)^{\otimes d}$.
Since \be \left(V^{\otimes d}\otimes (V^*)^{\otimes
d}\right)^{\rm{U}(\fg_{\bar 0})} \supset\left(V^{\otimes d}\otimes
(V^*)^{\otimes d}\right)^{\U}, \ee
 we have
\bea \left(V^{\otimes d}\otimes (V^*)^{\otimes d}\right)^{\U}
&=&j^{-1}\circ\rho_d(\C\fS_d). \eea This result may be stated more
explicitly as follows.
\begin{theorem} \cite{Sergeev01} \label{FFT1}
The vector space $\left(V^{\otimes d}\otimes (V^*)^{\otimes d}
\right)^{\U}$ is spanned by the following elements: \bea \sum_{a_1,
..., a_d} sgn(\sigma, a_1, \dots, a_d) v_{a_\sigma(1)}\otimes
v_{a_{\sigma(2)}}\otimes
...\otimes v_{a_{\sigma(d)}} \nonumber \\
\otimes {\bar v}_{a_d}\otimes {\bar v}_{a_{d-1}}\otimes ... {\bar
v}_{a_1}, \quad \forall \sigma\in\fS_d, \eea where $sgn(\sigma,
a_1, \dots, a_d)$ is a sign factor which is determined by the
restriction of $\sigma$ on the subset of odd indices in $\{a_1,
\cdots, a_d\}$ in such a way that if the restriction is even then
$sgn(\sigma, a_1, \dots, a_d)$ is $1$ and $-1$ otherwise.
\end{theorem}
We shall refer to both Theorems \ref{FFT} and  \ref{FFT1}
as the first fundamental theorem of the invariant theory of
the general linear supergroup.

\section{Superalgebras of Functions on the General Linear
Supergroup}\label{functions} \label{function-algebra}  We examine
properties of the Hopf superalgebra of regular functions on the
general linear supergroup in this section. The material presented here
is of critical importance for setting up the framework for studying
spherical functions. Some of the material can be extracted from
references \cite{SZ01, SZ02}.

Let $\U^0:=\{f\in\U^* \,|\, \ker f\
\mbox{contains a cofinite $\Z_2$-graded ideal of} \ \U\}$
be the finite dual \cite{Mo} of the universal enveloping
algebra $\U$ of $\g$. Standard Hopf algebra theory
\cite{MM, Mo} asserts that the Hopf
superalgebra structure of $\U$ induces a Hopf superalgebra structure on
$\U^0$. Denote by $m_{\circ}$, $\Delta_{\circ}$, $\epsilon_{\circ}$, and
$S_{\circ}$ the multiplication, comultiplication, counit, and antipode of
$\U^0$, respectively. The maps are defined for all
$f,\, g\in \U^0$ and $a,\, b\in \U$, by
\be
\langle m_{\circ}(f\otimes  g),\, a\rangle &=&
                            \langle f\otimes g, \, \Delta(a)\rangle \,, \\
\langle \Delta_{\circ}(f),\, a\otimes b\rangle
                                         &=&\langle f,\, a b \rangle \,, \\
\langle S_{\circ}(f),\,  a\rangle &=& \langle f,\, S(a)\rangle \,,
\ee
and $\1_{\U^0}=\epsilon$\,, $\epsilon_{\circ}=\1_{\U}$. Because $\U$ is
supercocommutative, $\U^0$ is  supercommutative. Recall that $S^2 = \id$
and hence also $S_{\circ}^2 = \id$. For convenience, we shall drop the
subscript $\circ$ from the notations for the multiplication,
comultiplication, and antipode of $\U^0$.

Let $\pi$ be a $\U$-representation of dimension $d<\infty$. Now for any $x\in\U$,
$\pi(x)$ is a $d\times d$-matrix.  We define a set of elements $\pi_{i j}\in\U^*$,
$i, j =1, 2, \ldots , d$, by
\be      \pi(x)&=& \left( \pi_{i j}(x)\right)_{i, j=1}^{d},   \quad
\forall x\in \U.\ee
The $\pi_{i j}$ will be called the matrix elements of $\pi$. It is easy to see
that the matrix elements of every finite-dimensional representation of $\U$ belong
to $\U^0$. Conversely, $\U^0$ is spanned by the matrix elements of all the
finite-dimensional representations of $\U$.
To see this,
we only need to consider an arbitrary non-zero element $f\in\U^0$.
Let $Ker$ be a graded cofinite ideal of $\U$ contained
in the kernel of $f$. Then $\U/Ker$ forms a left $\U$-module,
\be \U\otimes \U/Ker &\rightarrow&  \U/Ker, \\
       y\otimes (x+Ker) &\mapsto& y x+Ker.
\ee
Let $\{x_i + Ker\}$ be a basis of $\U/Ker$, and
denote by $f_{i j}$ the matrix elements of the associated representation
relative to this basis. Choose a set of complex numbers $c_i\in\C$ such that
$\1_{\U}+Ker=\sum_i c_i x_i +Ker$, where the set $\1_{\U}+Ker$ is not
contained in the kernel of $f$ since $f\ne 0$.
Then $f=\sum_{i, j} c_i \langle f, \, x_j\rangle  f_{j i}$.

We denote by  $t$  the $\U$-representation associated with the
contravariant vector module $V=L_{\epsilon_1}$ in the standard basis, and
denote its matrix elements by $t_{a b}\in \U^0$, $a, b\in{\bf I}$\,,
where $t_{a b}$ is even if $[a]+[b] = \bar{0}$, and odd otherwise.
Note that \be t_{a b}(E_{c d})&=& \delta_{a c}\delta_{b d}. \ee
Denote by $\bar{t}$ the covariant vector representation of $\U$
relative to the basis $\{{\bar v}_a \,|\, a\in{\bf I}\}$.
Let ${\bar t}_{a b}\in\U^0$\,, $a, b\in{\bf I}$, be the matrix
elements of ${\bar t}$. Then
\be {\bar t}_{a b} (E_{c d}) &=& - (-1)^{[a][b]+[b]}\delta_{b c} \delta_{a d}. \ee
Note that ${\bar t}_{a b}$ is even if
$[a]+[b] = \bar{0}$, and odd otherwise.

\begin{definition}\cite{SZ01}\label{function}
Let $\CG$ be the sub-superalgebra of $\U^0$ generated by
$\{ t_{a b}, \, \bar{t}_{a b}$ $|\, a, b\in{\bf I}\}$.
\end{definition}

The following relations hold in $\CG$
\begin{equation}
\sum_{c} t_{a c} {\bar t}_{b c} (-1)^{[c][a]+[b]} \,=\,\delta_{a b} \,,
                \quad
\sum_{c} {\bar t}_{c a} t_{c b} (-1)^{[b][c]+[c]} \,=\,\delta_{a b} \,,
\label{relation}
\end{equation}
because $t$ and ${\bar t}$ are dual representations of $\U$. More precisely,
the first relation states that the canonical tensor
$\sum_c v_c \otimes \bar{v}_c \in V \otimes V^*$ is $\U$-invariant, while
the second relation means that the dual pairing
$\langle \; , \: \rangle : V^* \otimes V \rightarrow \C$ is a $\U$-module
homomorphism.

$\CG$ has a bi-superalgebra structure, with
the co-multiplication defined by
\be \Delta( {t}_{a b} ) \,=\, \sum_{c\in{\bf I}}
      (-1)^{([c]-[a])([c]-[b])} {t}_{a c}\otimes {t}_{c b},\\
 \Delta( \bar{t}_{a b} ) \,=\, \sum_{c\in{\bf I}}
      (-1)^{([c]-[a])([c]-[b])} \bar{t}_{a c}\otimes \bar{t}_{c b}.
      \ee
Let us also denote by $S$ the antipode of $\U^0$. By using the definition
of dual modules we can show that
\bea
   S(t_{a b}) = (-1)^{[a][b]+[a]} {\bar t}_{b a},
                &\quad&
   S({\bar t}_{a b}) = (-1)^{[a][b]+[b]} t_{b a}. \label{antipode}
\eea
The following result was proven in \cite{SZ01}.
\begin{proposition}\cite{SZ01}
(1). $\CG$ forms a Hopf sub-superalgebra of $\U^0$.\\
(2). $\CG$ is dense in $\U^*$ in the following sense:
for every non-zero element $x\in \U$, there exists some $f\in \CG$ such that
$\langle f,\, x\rangle \ne 0$.
\end{proposition}

Let $\Lambda$ denote a finite dimensional Grassmann algebra. Recall
that the general linear supergroup $GL(m|n, \Lambda)$ over $\Lambda$
is the group of even invertible $(m+n)\times(m+n)$-matrices with
entries in $\Lambda$. It was shown in \cite{SZ02} that $GL(m|n,
\Lambda)$ can be reconstructed from $\CG$ in the following way. The
$\Z_2$-graded vector space $\Hom_{\C}\left(\CG, \Lambda\right)$ has a
natural superalgebra structure, with the multiplication defined for
any $\phi$ and $\psi$  by \bea (\phi\psi)(f) := \sum_{(f)}
(-1)^{[f_{(1)}][\psi]} \phi(f_{(1)}) \psi(f_{(2)}), \quad \forall
f\in\CG, \label{product} \eea where we have used Sweedler's notation
expressing the co-multiplication $\Delta(f)$ of any $f\in \CG$ by
$\sum_{(f)} f_{(1)} \otimes f_{(2)}$.
\begin{theorem} \label{supergroup} \cite{SZ02}
Let $G_\C:=\{ \mbox{superalgebra homomorphisms}\ \CG\rightarrow
\Lambda\}$. Then with the multiplication defined by (\ref{product}),
the set $G_\C$ forms a group, which is isomorphic to $GL(m|n,
\Lambda)$.
\end{theorem}
We shall not repeat the proof of the Theorem here, but merely point
out that the inverse $\alpha^{-1}$ of any element $\alpha\in G_\C$ is
given by $\alpha^{-1}(f)=\alpha(S(f))$, for all $f\in\CG$.

We shall refer to the elements of $\CG$ as the regular functions on
the general linear supergroup. We now consider their properties. Note
that there exists two natural left actions $dR$ and $dL$ of $\U$ on
$\CG$ respectively corresponding to the left and right translations.
For all $x\in\U$, $f\in\CG$, \bea dR_x(f)&=&\sum_{(f)} (-1)^{[x][f]} \,
         f_{(1)} \, \langle f_{(2)},\, x\rangle \,,\nonumber \\
dL_x(f)&=&\sum_{(f)} (-1)^{[x]} \, \langle f_{(1)},\, S(x)
                                           \rangle f_{(2)}. \label{circ}
\eea
Equivalently, the equations in
(\ref{circ}) can be rewritten in the form
\begin{equation}
 \langle dR_x(f), \, y\rangle  = (-1)^{[x]([f]+[y])} \langle f, yx \rangle \,,\quad
 \langle dL_x(f), \, y\rangle = (-1)^{[x][f]} \langle f, S(x)y \rangle \,, \label{circalt}
\end{equation}
for all $x,y \in \U$ and $f \in \CG$.
Straightforward calculations show that each of $d L$ and $d R$ indeed
converts $\CG$ into a (graded) left $\U$-module.
With respect to this module structure the product map of $\CG$ is
a $\U$-module homomorphism, and the unit element of $\CG$ is
$\U$-invariant.  Take $d L$ as an example, we have
\bea  \sum_{(x)} (-1)^{[x_{(2)}] [f]} d L_{x_{(1)}}(f)\otimes
d L_{x_{(1)}}(g) &\mapsto& d L_x(f g),
\quad \forall f, \, g\in \U^0, \ x\in\U, \nonumber\\
d L_x (\epsilon) &=&  \epsilon(x) \epsilon, \quad \forall x\in\U.
\label{modalg}\eea
This is saying that each of the actions $d L$ and $d R$ converts $\CG$
into a left $\U$-module superalgebra \cite{Mo}.
The two actions supercommute as can be easily checked. Thus
$\CG$ forms a left $\U\otimes \U$-module algebra, with the action
\be (x\otimes y) f &=& d L_x d R_y (f), \quad \forall x, y\in \U, \
f\in\CG.\ee
The fact that the product map in
$\CG$ is a module homomorphism means that the operators $dR_x$ and $dL_x$
behave as some sort of generalized superderivations. In particular, if
$x \in \g$, they are superderivations.

To better understand the structure of $\CG$, we let
$X= V\otimes V^*$ and $\bar{X}=V^*\otimes V$. Using the
standard bases of $V$ and $V^*$ we manufacture the bases
$\{x_{a b}:= v_b\otimes {\bar v}_a\}$ and
$\{{\bar x}_{a b}:={\bar v}_b \otimes v_a\}$ for
$X$ and $\bar X$ respectively. Denote by $\CXbX$ the $\Z_2$-graded
symmetric algebra of $X\oplus {\bar X}$. Then
$\CXbX$ as an associative superalgebra can be described more explicitly
as  generated by $x_{a b}$\,,
${\bar x}_{a b}$\,, $a, b\in{\bf I}$\,, subject to the relations
\be
x_{a b} x_{c d} &=& (-1)^{([b]-[a])([d]-[c])} x_{c d} x_{a b}\,,\\
x_{a b} {\bar x}_{c d} &=&
                    (-1)^{([b]-[a])([c]-[d])} {\bar x}_{c d} x_{a b}\,,\\
{\bar x}_{a b} {\bar x}_{c d} &=&
                    (-1)^{([a]-[b])([c]-[d])} {\bar x}_{c d} {\bar x}_{a b}\,.
\ee
The generators $x_{a b}$ and ${\bar x}_{a b}$ are even if
$[a]+[b] = \bar{0}$, and odd otherwise. Stated differently,
the $2(m^2+n^2)$ even generators generate a polynomial algebra,
the $4m n$ odd generators generate a Grassmann algebra with the standard
grading, and $\CXbX$ is the tensor product of the two.
Let $\cJ$ be the (graded) ideal of $\CXbX$ generated by the following elements:
\bea
\sum_{c} x_{a c} {\bar x}_{b c} (-1)^{[c][a]+[b]} - \delta_{a b} \,,
                &\quad&
\sum_{c} {\bar x}_{c a} x_{c b} (-1)^{[b][c]+[c]} - \delta_{a b} \,,
                \quad a, b\in{\bf I} \,.  \label{ideal}
\eea
We have the following theorem.

\begin{theorem} \cite{SZ01} \label{poly}
The assignments
$ x_{a b}\mapsto t_{a b} \,, \quad {\bar x}_{a b}\mapsto {\bar t}_{a b} \,,
                                                \quad a, b\in{\bf I} $
specify a well-defined superalgebra isomorphism
$\jmath:  \CXbX/{\cJ}\rightarrow\CG.$
\end{theorem}

Define two left $\U$-actions on $X\oplus \bar{X}$
\be \Phi: \U\otimes (X\oplus \bar{X}) \rightarrow X\oplus \bar{X},
&\quad& u\otimes w \mapsto \Phi(u) w, \\
 \Psi: \U\otimes (X\oplus \bar{X}) \rightarrow X\oplus \bar{X},
 &\quad& u\otimes w \mapsto \Psi(u) w.\ee
by \be \Phi(u)(v_b\otimes \bar{v}_a) &=& (-1)^{[u]} u v_b\otimes
\bar{v}_a, \\
\Psi(u)(v_b\otimes \bar{v}_a) &=& (-1)^{[u][b]}  v_b\otimes
u \bar{v}_a,\\
\Phi(u) (\bar{v}_b\otimes v_a) &=& (-1)^{[u]} u \bar{v}_b\otimes v_a,\\
\Psi(u) (\bar{v}_b\otimes v_a) &=& (-1)^{[u]([u]+[b])}  \bar{v}_b\otimes u v_a,
\quad u\in\U.
\ee
These actions super-commute, and can both be extended to left $\U$-actions on $\CXbX$ by
\be \Phi(x)(p_1 p_2) &=& \sum (-1)^{[x_{(2)}][p_1]}\left(\Phi(x_{(1)})p_1 \right)
 \left(\Phi(x_{(2)}) p_2 \right),\\
\Psi(x)(p_1 p_2) &=& \sum (-1)^{[x_{(2)}][p_1]}\left(\Psi(x_{(1)})p_1 \right)
 \left(\Psi(x_{(2)}) p_2 \right),
\ee
where $p_1, p_2\in \CXbX$ and $x\in\U$. This gives rise to a $\U\otimes\U$-module
algebra structure on $\CXbX$. Note that the
$\U\otimes\U$-action leaves the ideal $\cJ$ invariant.
Thus we have the following proposition.
\begin{proposition}\label{quotient}
The map $\jmath:  \CXbX/{\cJ} \rightarrow \CG$ of Theorem \ref{poly} is a
$\U\otimes\U$-module algebra isomorphism, with
\be \jmath\left((\Psi(x)\otimes\Phi(y))p\right) &=& (d L_x\otimes d R_y)\jmath(p),
\quad \forall x, y\in\U, \ p\in \CXbX.  \ee
\end{proposition}

\section{Homogeneous Superspaces and Spherical
Functions}\label{sect-space} Recall the following well known fact
in the context of classical homogeneous spaces: if $H$ is a
compact semi-simple Lie group, and $H_\C$ is its complexification,
then for any parabolic subgroup $Q$ of $H_\C$, we have $H_\C/Q =
H/R$, where $R$ is the intersection of the Levi factor of $Q$ with
$H$. We shall imitate this construction in the algebraic setting
for Lie supergroups. For this we need to discuss real forms of the
complex general linear superalgebra and the general linear
supergroup.

\subsection{Real Forms}\label{real}
Let us begin by briefly discussing the notion of Hopf
$\ast$-superalgebras \cite{Zh98}. A $\ast$-superalgebraic
structure on an associative superalgebra $A$ is a conjugate linear
anti-involution $\theta: A\rightarrow A$: for all $x, y\in A$,
$c,c'\in\C$, \be \theta(c x + c' y)={\bar c}\theta(x) + {\bar c'}
\theta(y),\quad& \theta(x y)=\theta(y)\theta(x), \quad&
\theta^2(x)=x. \ee Note that the second equation does not involve
any sign factors as one would normally expect of superalgebras. We
shall sometimes use the notation $(A, \theta)$ for the
$\ast$-superalgebra $A$ with the $\ast$-structure $\theta$.  Let
$(B, \theta_1)$ be another associative $\ast$-superalgebra. Now
$A\otimes B$ has a natural superalgebra structure, with the
multiplication defined for any $a, a'\in A$ and $b, b'\in B$ by
\be (a\otimes b) (a'\otimes b') &=& (-1)^{[b][a']}a a'\otimes b
b', \ee where $(-1)^{[b][a']}$ is the usual sign factor required
for exchanging positions of odd elements. Furthermore, the
following conjugate linear map \bea \theta\star\theta_1: a\otimes
b \mapsto (1\otimes \theta_1(b)) (\theta(a)\otimes 1) =
(-1)^{[a][b]} \theta(a)\otimes\theta_1(b) \label{star}\eea defines
a $\ast$-superalgebraic structure on $A\otimes B$.

Let us assume that $A$ is a Hopf superalgebra with co-multiplication
$\Delta$, co-unit $\epsilon$ and antipode $S$.  If the
$\ast$-superalgebraic structure $\theta$ satisfies  \be
(\theta\star\theta)\Delta=\Delta\theta, &\quad&
\theta\epsilon=\epsilon\theta,  \ee then $A$ is called a Hopf
$\ast$-superalgebra.  Now $$\sigma:=S \theta$$ satisfies
$\sigma^2=id_A$, as follows from the definition of the antipode.

Let $A^0$ denote the finite dual of $A$, which has a natural Hopf
superalgebraic structure. We shall still use $\Delta$ and $S$ to
respectively denote the co-multiplication and antipode of $A^0$, but
write its co-unit as $\epsilon_o$. If $A$ is a Hopf $\ast$-superalgebra
with the Hopf $\ast$-superalgebraic structure $\theta$, then $\sigma=S
\theta$ induces a map $\omega: A^0\rightarrow A^0$ defined for any
$f\in A^0$ by \bea \langle \omega(f), x\rangle = \overline{\langle f,
\sigma(x)\rangle }, \quad \forall x\in A. \label{ast}\eea
\begin{lemma}\label{dual-ast}
The map $\omega: A^0\rightarrow A^0$ defined by \eqref{ast} gives rise
to a Hopf $\ast$-superalgebraic structure on $A^0$.
\end{lemma}
\begin{proof}
It is clear that $\omega$ is conjugate linear. Also, $\sigma^2=id_A$
implies $\omega^2=id_{A^0}$. For all $f, g\in A^0$, $x, y\in A$, we
have \be \langle \omega(f g), x\rangle &=& \overline{\langle f g,
\sigma(x)\rangle }
=\overline{\langle f\otimes g,  (S\otimes S)(\theta\star\theta)\Delta'(x)\rangle }\\
&=&  (-1)^{[f][g]}\langle \omega(f)\otimes\omega(g),  \Delta'(x)
\rangle = \langle \omega(g)\omega(f), x\rangle,  \ee that is,
$\omega(f g)=\omega(g)\omega(f)$. Define $\omega\star\omega$ as in
\eqref{star}, we have \be \langle (\omega\star\omega)\Delta(f),
x\otimes y \rangle &=&(-1)^{[x][y]} \overline{\langle \Delta(f) ,
\sigma(x)\otimes\sigma(y)\rangle } = \overline{\langle f , \sigma(x y)\rangle }\\
&=&\langle \omega(f) , x y\rangle=\langle \Delta\omega(f), x\otimes y
\rangle, \ee that is $(\omega\star\omega)\Delta(f)=\Delta\omega(f).$
It is easy to show that $\omega$ also satisfies all the other minor
requirements to qualify as a Hopf $\ast$-superalgebraic structure on
$A^0$.
\end{proof}

The universal enveloping algebra of the general linear
superalgebra admits many Hopf $\ast$-superalgebraic structures.
Let us fix one Hopf $\ast$-superalgebraic structure $\theta:
\rm{U}(\fg)\rightarrow \rm{U}(\fg)$ here. As $\fg$ is canonically
embedded in $\rm{U}(\fg)$, the restriction of $\theta$ to $\fg$
defines a conjugate anti-involution of the Lie superalgebra. Let
$\fg_{\bar 0}^\sigma$ and $\fg_{\bar 1}^\sigma$ be the fixed point
sets of $\fg_{\bar 0}$ and $\fg_{\bar 1}$ under $\sigma$
respectively. Let $\fgu\subset\fg$ be the real span of $\fg_{\bar
0}^\sigma\cup \sqrt{i} \fg_{\bar 1}^\sigma$. Then $\fgu$ forms a
real Lie superalgebra, which is a real form of $\fg$. However,
note that the $\sigma$-invariants of $\fg$ do not form a real
subalgebra of $\fg$ if $\fg_{\bar 1}^\sigma$ is non-trivial. This
is the reason for us to consider $\fgu$ instead.

Denote by $\rm{U}^\R(\fgu)$ the real universal enveloping algebra
of $\fgu$, which is embedded in $\U$ in the natural way.
Furthermore,
$$\rm{U}(\fg)= \C\otimes_{\R}\rm{U}^\R(\fgu).$$
By Lemma \ref{dual-ast}, the Hopf $\ast$-superalgebraic structure
$\theta$ induces a Hopf $\ast$-superalgebraic structure $\omega:
\CG\rightarrow\CG$ on $\CG$.  By using the embedding of the real
associate superalgebra $\rm{U}^\R(\fgu)$ in $\rm{U}(\fg)$, we can
see that $f\in\CG$ vanishes if and only if $\langle f, x\rangle
=0$, for all $x\in \rm{U}^\R(\fgu)$. Therefore elements of $\CG$
can be considered as complex valued functionals on the real
superalgebra $\rm{U}^\R(\fgu)$. From this point of view, we should
interpret $\CG$ as the $\ast$-superalgebra of functions on some
real supergroup $G$. Now let us make this discussion more precise.

Let $\Lambda$ be the complex Grassmann algebra introduced in section
\ref{function-algebra}. Let $ ^-: \Lambda\rightarrow \Lambda$ be a
`complex conjugation operation' on supernumbers (i.e., $(\Lambda,\
^-)$ is a $\ast$-superalgebra). Theorem \ref{supergroup} shows that
all the superalgebra homomorphisms $\CG\rightarrow \Lambda$ form a
supergroup $G_\C$, which is isomorphic to $GL(m|n, \Lambda)$. A
homomorphism $\alpha: \CG\rightarrow \Lambda$ will be called a
$\ast$-superalgebra homomorphism if it preserves the
$\ast$-superalgebraic structures in the sense that
$\alpha(\omega(f))=\overline{\alpha(f)}$, for all $f\in\CG$. The
following result can be easily proven.
\begin{lemma}
If an element $\alpha$ of  $G_\C$ is a $\ast$-superalgebra
homomorphism, then its inverse is also a $\ast$-superalgebra
homomorphism. The product of any two $\ast$-superalgebra homomorphisms
in $G_\C$ is again a $\ast$-superalgebra homomorphism.
\end{lemma}
\begin{proof}
We shall prove the first statement only. The second one can be shown
in a similar way. Recall that the inverse of $\alpha\in G_\C$ is
defined by \be \langle\alpha^{-1}, \ f \rangle &=& \langle\alpha, \
S(f) \rangle, \quad \forall f\in\CG. \ee  Now if $\alpha$ is a
$\ast$-superalgebra homomorphism, then for all $f\in \CG$, \be
\langle\alpha^{-1}, \ \omega(f) \rangle = \langle\alpha, \ S\omega(f)
\rangle = \overline{\langle\alpha, \ S(f)
\rangle}=\overline{\langle\alpha^{-1}, \ f \rangle}.\ee This shows
that $\alpha^{-1}$ is indeed a $\ast$-superalgebra homomorphism.
\end{proof}
Introduce the map $\check\theta: G_\C\rightarrow G_\C$ defined by \be
\langle \check\theta(\alpha), \ f \rangle &=& \overline{\langle\alpha,
\ \omega S(f) \rangle}, \quad \forall f\in\CG. \ee  We need to show
that the image of $\check\theta$ indeed lies in $G_\C$. For any $f,
g\in\CG$, we have \be \langle \check\theta(\alpha), \ f g \rangle &=&
(-1)^{[f][g]} \overline{\langle\alpha, \ \omega S(f) \omega S(g)
\rangle}\\ &=& (-1)^{[f][g]} \overline{\langle\alpha\otimes\alpha,
\ \omega S(f)\otimes \omega S(g) \rangle}\\
&=& (-1)^{[f][g]} \overline{\langle\alpha, \ \omega S(f) \rangle
\langle\alpha, \ \omega S(g) \rangle}\\
&=& \overline{\langle\alpha, \ \omega S(g) \rangle
\langle\alpha, \ \omega S(f)\rangle  }\\
&=& \overline{\langle\alpha, \ \omega S(f)\rangle } \cdot
\overline{\langle\alpha, \ \omega S(g) \rangle}\\
&=& \langle \check\theta(\alpha), \ f\rangle
\langle\check\theta(\alpha), \ g \rangle.  \ee Therefore,
$\check\theta(\alpha)$ is a superalgebra homomorphism from $\CG$ to
$\Lambda$, thus is indeed an element of $G_\C$.
\begin{definition}
$G:=\{\ast-\rm{superalgebra \ homomorphism} \ \CG\rightarrow \Lambda
\}$.
\end{definition}
\begin{theorem} $G$ forms a subgroup of $G_\C$. Furthermore,
$\check\theta(\alpha)=\alpha^{-1}$ for all $\alpha\in G$.
\end{theorem}
\begin{proof}
The fact that $G$ forms a subgroup immediately follows from the above
lemma. If $\alpha\in G$, we have \be \langle \check\theta(\alpha), \ f
\rangle &=& \overline{\langle\alpha, \ \omega S(f) \rangle} =
\langle\alpha, \ S(f) \rangle\\
&=& \langle\alpha^{-1}, \ S(f) \rangle, \quad \forall f\in\CG.\ee This
confirms the second claim.
\end{proof}

\subsection{Spherical functions on homogeneous
superspaces}\label{subsect-space} Hereafter we fix a Hopf
$\ast$-superalgebraic structure $\theta$ for $\U$, which is
defined for all the generators by \bea \theta: E_{a b} \mapsto
E_{b a}. \label{compact-gl} \eea The associated real form
${\mathfrak{gl}}(m|n;\C)^{\sigma, \sqrt{i}}$ of the general linear
superalgebra is one of the `compact' real forms of the general
linear superalgebra, which probably deserves the notation
$\fu(m|n)$ because it contains the maximal even subalgebra
$\fu(m)\oplus\fu(n)$. (The unitarizable representations of this
compact real form comprise of the tensor powers of the natural
representation, while the unitarizable representations of the
other compact real form are the duals of these representations
\cite{Zh89}.) Direct calculations can show that the Hopf
$\ast$-superalgebraic structure on $\CG$ induced by $\theta$ is
given by \bea \omega(t_{a b})= (-1)^{[b]([a]+[b])} {\bar t}_{a b},
\quad \omega(\bar t_{a b})= (-1)^{[b]([a]+[b])} {t}_{a b}.
\label{compact-G}\eea  The real supergroup $G$ has body
$U(m)\times U(n)$.

Let $\fp$ be a parabolic subalgebra of $\fg$ with Levi factor
$\fl$.  Let $\fk = \fl^{\sigma, \sqrt{i}}$ be the real form of
$\fl$, which is a subalgebra of $\fgu$. Denote by $\rm{U}^\R(\fk)$
the universal enveloping algebra of $\fk$ over $\R$. Note that
$\rm{U}^\R(\fgu)$ inherits a real Hopf superalgebra structure from
$\rm{U}(\fg)$, and $\rm{U}^\R(\fk)$ inherits a real Hopf
superalgebra structure from $\rm{U}^\R(\fgu)$. Let us introduce
the following definition.
\begin{definition}
\bea
  \CKG &:=& \left\{ f\in\CG \,|\, dL_k(f) =\epsilon(k) f,
 \ \forall\, k\in\URk \right\}.  \label{space}
\eea\end{definition} Note the following obvious fact, which will
be used immediately below: \bea
  \CKG &:=& \left\{ f\in\CG \,|\, dL_k(f) =\epsilon(k) f,
 \ \forall\, k\in\Ul \right\}.  \label{complex-setting-a}
\eea We have the following lemma.
\begin{lemma}\label{subalgebra} $\CKG$ forms a $\ast$-subalgebra of $\CG$,
which is also a left co-ideal of $\CG$.
\end{lemma}
\begin{proof}
Since $\Ul$ is a Hopf subalgebra of $\rm{U}(\fg)$, we have
$\Delta(k)=\sum_{(k)}k_{(1)}\otimes k_{(2)}$ $\in\Ul\otimes\Ul$
for all $k\in\Ul$. If $a, b \in\CKG$, then by
(\ref{complex-setting-a}), \be d L_k ( a b )&=&\sum
(-1)^{[a_{(2)}][b_{(1)}]+[k]}
\langle a_{(1)} b_{(1)},\ S(k)\rangle a_{(2)} b_{(2)}\\
&=&\sum (-1)^{[a_{(2)}][k_{(1)}]+[k]}
\langle a_{(1)}, S(k_{(2)})\rangle \langle b_{(1)},\ S(k_{(1)})\rangle
a_{(2)} b_{(2)}\\
&=&\sum (-1)^{[k]} \epsilon(k_{(1)}) \langle a_{(1)},
S(k_{(2)})\rangle a_{(2)} b =\epsilon(k) a b,  \quad \forall\,
k\in\Ul. \ee Thus $a b\in\CKG$.

Given any $f\in\CKG$, we have $dL_k(\omega(f))
=\omega(dL_{\theta(k)}(f))(-1)^{[k]([k]+[f])}$ for all $k\in\Ul$.
As $\Ul$ is $\theta$ invariant, we have
$dL_{\theta(k)}(f)=\overline{\epsilon(k)} f$. Thus
$$dL_k(\omega(f)) =\epsilon(k)\omega(f), \quad \forall k\in\Ul.$$
Also, a straightforward calculation shows that
$$(dL_k\otimes \id)\Delta(f)=\epsilon(k)\Delta(f), \quad \forall k\in\Ul.$$
Thus $\CKG$ is a left co-ideal. This completes the proof.
\end{proof}

The subalgebra $\CKG$ consists of the elements of $\CG$ which are
invariant with respect to $\URk$ under `left translation'.
Following the general philosophy of non-commutative geometry
\cite{Co}, we may take the viewpoint that $\CKG$ defines an
algebraic homogeneous superspace \cite{Ma88}. We shall refer to
$\CKG$ as the superalgebra of functions on the homogeneous
superspace. Also a word about the notation $\CKG$: here $K$ is
used to indicate some real sub-supergroup of $G$ with Lie
superalgebra $\fk$.

\begin{remark}
Since $\CG$ and $\CKG$ are all $\ast$-superalgebras, their elements
are in general not `holomorphic functions'  on the supergroup. This is
a particularly welcome fact, as it indicates that our construction can
lead to analogues of compact complex super manifolds like projective
superspaces. As is well known from the Gelfand-Naimark theorem, the
continuous functions on a compact manifold determine the manifold
completely, even when the manifold is complex, where all the
holomorphic functions are constants.
\end{remark}

\begin{remark}
In the quantum group context, one usually considers left or right
co-ideal subalgebras of the algebra of functions \cite{KD, No, DN,
Le} in the place of $\CKG$.  By Lemma \ref{subalgebra} $\CKG$
forms a left co-ideal subalgebra of $\CG$.
\end{remark}

Because the two left actions $d R$ and $d L$ of $\U$ on $\CG$
super-commute, the subalgebra $\CKG$ of $\CG$ forms a left module
algebra over $\U$ under the action $d R$. We shall study the $d
R(\URk)$-invariant subspace of $\CKG$. Let us first generalize the
definition of zonal spherical functions \cite{Ta} to the
supergroup setting. We shall refer to elements of the following
space as spherical functions on the homogeneous superspace.
\begin{definition}
\bea \CKGK&:=&\{ f\in\CKG \mid d R_k (f) =\epsilon(k) f,
 \ \forall\, k\in\URk \} \label{spherical}
\eea \end{definition}
Similar arguments as those in the proof of Lemma
\ref{subalgebra} show that
\begin{lemma}
The subspace $\CKGK$ forms a $\ast$-subalgebra of $\CKG$.
\end{lemma}
Obviously  \bea
 \CKGK&=&\{ f\in\CKG \mid d R_x (f) =\epsilon(x) f,
 \ \forall\, x\in \rm{U}(\fl) \}, \label{complex-setting}
\eea where $\fl$ is the complexification of $\fk$. The fact will
be used in the next subsection to prove Theorem
\ref{construction}.

\subsection{Structure of superalgebra of spherical functions}\label{main}
Let $\fl$ be a reductive subalgebra of $\fg$ generated by $E_{a
a}$, $a\in {\bf I}$, and $E_{c,c+1}, \ E_{c+1, c}$ with $c$
belonging to some proper subset of ${\bf I}\backslash\{m+n\}$. As
in the last subsection, we let $$\fk=\fl^{\sigma, \sqrt{i}}.$$ See
Remark \ref{homogeneous-space} for further discussions on this
choice of $\fk$. The main result here is Theorem
\ref{construction}, which enables us to obtain the superalgebras
$\CKG$ and $\CKGK$ from the invariants of $\CXbX$. An explicit
description of the generators of these superalgebras will also be
given in Lemmas \ref{gCKG} and \ref{gCKGK}.
\begin{theorem} \label{construction} When $\fk=\fl^{\sigma, \sqrt{i}}$, we have
\bea \CKG&=& \left\{ \jmath(p) | p\in
\CXbX^{\Psi(\Ul)}\right\}, \quad {\text and}\nonumber\\
\CKGK&=& \left\{ \jmath(p) | p\in
\CXbX^{\Psi(\Ul)\otimes\Phi(\Ul)}\right\}. \eea
\end{theorem}
The remainder of this subsection is devoted to the proof of
Theorem \ref{construction}. The proof is carried out in two steps.
We first show that the theorem holds when $\fl=\fk_\C$ is even,
that is, when $\fl$ is a reductive Lie subalgebra of $\fg$. Then
we use this fact to prove the general case. In the process of
proving the theorem, we also establish Lemmas \ref{gCKG} and
\ref{gCKGK}. We mention that equations  (\ref{complex-setting-a})
and (\ref{complex-setting}) will be used repeatedly in the proof
without further warning.

\subsubsection{Proof of Theorem \ref{construction} for $\fl$ even}
In this case we can find a set of positive integers $k_i$, $i=1, 2,
..., r, r+1,  ..., s$ such that $\sum_{i=1}^r k_i =m$, $\sum_{j=r+1}^s
k_j=n$, and $\fl= \oplus_{i=1}^s \mathfrak{gl}(k_i).$ More explicitly,
\be \fl&=&\left\{ \left.\left(\begin{array}{r c l}
A_1 & & 0  \\
    &\ddots& \\
0  & & A_s \end{array} \right)\right|
A_i\in\mathfrak{gl}(k_i)\right\}\subset \fg.\ee Proposition
\ref{quotient} implies the following short exact sequence \be
0\longrightarrow \cJ \longrightarrow
\CXbX\stackrel{\jmath}{\longrightarrow} \CG\longrightarrow 0\ee in the
category of $\Ul\otimes\Ul$-module superalgebras. Since the various
$\Ul$ and $\Ul\otimes \Ul$ actions on $\cJ$, $\CXbX$ and $\CG$ are all
semi-simple, we have the following short exact sequences of
$\Ul\otimes\Ul$-modules \be &0\longrightarrow \cJ^{\Psi(\rm{U}(\fl))}
\longrightarrow \CXbX^{\Psi(\rm{U}(\fl))}{\longrightarrow}
\CG^{dL_{\rm{U}(\fl)}}\longrightarrow 0, &\\
&0\longrightarrow \cJ^{\Psi(\rm{U}(\fl))\otimes\Phi(\rm{U}(\fl))}
\longrightarrow
\CXbX^{\Psi(\rm{U}(\fl))\otimes\Phi(\rm{U}(\fl))}{\longrightarrow}
\CG^{d L_{\rm{U}(\fl)}\otimes d R_{\rm{U}(\fl)}}\longrightarrow 0,&\ee
where $\CKG=\CG^{d L_{\rm{U}(\fl)}}$ and $\CKGK=\CG^{d
L_{\rm{U}(\fl)}\otimes d R_{\rm{U}(\fl)}}$. These are also short exact
sequences of $\rm{U}(\fl)\otimes\rm{U}(\fl)$-module algebras, thus
they imply the claims of Theorem \ref{construction} in the case under
consideration.

Let us now describe the algebras $\CKG$ and $\CKGK$ more carefully.
Set $l_i=\sum_{t=1}^i k_t$. Recall that $\CXbX$ is the symmetric
algebra in $X\oplus {\bar X}$ where $X= V\otimes {\bar V}$ and ${\bar
X}={\bar V} \otimes V$. Restricted to a $\rm{U}(\fl)$-module, $V$
decomposes into \be V&=& \oplus_{i=1}^s V_i^{(k_i)}.\ee The ideal
$\mathfrak{gl}(k_i)$ of $\fl$ acts on $V_i^{(k_i)}$ by the natural
representation,  and acts  on all other submodules trivially. There is
also an analogous decomposition of the restriction of $\bar V$ to a
$\rm{U}(\fl)$-module. By applying the first fundamental theorem of the
invariant theory of the general linear group \cite{Ho}, we obtain that
the subalgebra $\CXbX^{\Psi(\rm{U}(\fl))}$ of $\CXbX$ is generated by
\be {\hat C}^{(i)}_{a b}&:=& \sum_{c=1+l_{i-1}}^{l_i} x_{c a} {\bar
x}_{c b}, \quad  i=1, 2, ..., s, \quad a, b\in \bf{I}.\ee It then
immediately follows that $\CKG$ is generated by \be C^{(i)}_{a
b}&:=\jmath({\hat C}^{(i)}_{a b})&= \sum_{c=1+l_{i-1}}^{l_i} t_{c a}
{\bar t}_{c b}, \quad i=1, 2, ..., s,\quad a, b \in \bf{I}.\ee

Note that the $C^{(i)}_{a b}$ are not algebraically independent,
for example, for $a, b \in {\bf I}$ the following hold
\bea \sum_{i=1}^s C_{a b}^{(i)}
(-1)^{[l_i]} = \delta_{a b}, &\quad& \sum_{a=1}^{m+n} C_{a
b}^{(i)} = k_i. \eea
Thus the elements of the set $\{ C_{a b}^{(i)} \mid i\ne s; a, b\in{\bf I}\}$
can also generate $\CKG$. By using the fact that $t_{ab}$
and ${\bar t}_{cd}$ super-commute, one can verify the following
proposition easily.

\begin{proposition} We have
\bea C_{ab}^{(i)}C_{cd}^{(j)}= (-1)^{([a]+[b])([c]+[d])} C_{cd}^{(j)}C_{ab}^{(i)}, \eea
in particular, if $[a]+[b]=1$ then $(C_{ab}^{(i)})^2 =0$.  Thus for fixed $i$,
there is an onto algebra homomorphism $\C [X] \rightarrow <C_{ab}^{(i)}|a,b\in \bf{I}>$,
where $\C [X]$ is the subalgebra of $\CXbX$ generated by $X$, and $<C_{ab}^{(i)}
| a, b \in \bf{I}>$ is the subalgebra of $\CKG$ generated by
$\{ C_{ab}^{(i)} | a, b \in \bf{I} \}$.
\end{proposition}

In a similar way we can show that $\CXbX^{\Psi(\rm{U}(\fl))\otimes
\Phi(\Ul)}$ is generated by \be {\hat C}^{(i, j)}&:=&
\sum_{a=1+l_{j-1}}^{l_j}\sum_{c=1+l_{i-1}}^{l_i} x_{c a} {\bar x}_{c
a}, \quad i, j=1, 2, ..., s,\ee  and $\CKGK$ is generated by \be
{C}^{(i, j)}&:=& \jmath({\hat C}^{(i, j)})
=\sum_{a=1+l_{j-1}}^{l_j}\sum_{c=1+l_{i-1}}^{l_i} t_{c a} {\bar t}_{c
a}, \quad i, j=1, 2, ..., s.\ee Again, the $C^{(i, j)}$ are not
algebraically independent, for example, \bea
 \sum_{i=1}^s C^{(i, j)} (-1)^{[l_i]} = k_j, &\quad& \sum_{j=1}^s C^{(i, j)}
(-1)^{[l_j]} = k_i.\eea Thus the elements of the set $\{{C}^{(i, j)}
\mid i, j\ne r\}$ generate $\CKGK$.

\subsubsection{Proof of Theorem \ref{construction} for generic $\fl$}
The most general form of $\fl$ is as follows. There exists a set of
positive integers $k_i$ as in the last subsection such that
$$\fl=\left(\oplus_{i=1}^{r-1}\mathfrak{gl}(k_i) \right)\oplus
\mathfrak{gl}(k_r|k_{r+1}) \oplus
\left(\oplus_{j=r+2}^{s}\mathfrak{gl}(k_j) \right).$$  More
explicitly, we have \be \fl&=&\left\{ \left.\left(\begin{array}{c c c c
c c c}
A_1&      &        &   &   &  & \\
   &\ddots&        &   &   & 0 & \\
   &      & A_{r-1}&   &   &  & \\
   &      &        & B &   &  &\\
& & & &A_{r+2}&      &          \\
& 0& & &      &\ddots&           \\
& & & &      &      & A_{s}      \\
\end{array}
\right)\right|\begin{array}{r c l}
A_i&\in&\mathfrak{gl}(k_i),\\
B&\in&\mathfrak{gl}(k_r|k_{r+1}) \end{array}\right\}.\ee Note that
$\fl$ contains the maximal even subalgebra
$\fl_0=\oplus_{i=1}^{s}\mathfrak{gl}(k_i).$

We first consider the subalgebra $\CG^{dL_{\rm{U}(\fl_0)}}$ of
$\CG$. By using results of the last subsection, we can immediately
see that $\CG^{d L_{\rm{U}(\fl_0)}}$ is generated by the elements
of $\{C_{a b}^{(i)} | i\ne r; \ a, b\in {\bf I}\}$. Now \be
\CKG&=& \{f\in \CG^{d L_{\rm{U}(\fl_0)}} | d L_{E_{m m+1}}(f)=d
L_{E_{m+1, m}}(f)=0 \}. \ee We shall show that $\CKG$ is generated
by $\{C_{a b}^{(i)} | i\ne r, r+1; \ a, b\in {\bf I}\}$.

Note that all the elements of $\{C_{a b}^{(i)} | i\ne r; \ a, b\in {\bf I}\}$
are annihilated by $d L_{E_{m m+1}}$ and $d L_{E_{m+1, m}}$ except
for $C_{a b}^{(r+1)}$, for which we have
\bea d L_{E_{m m+1}}(C_{a b}^{(r+1)})&=& -(-1)^{[a]+[b]} t_{m+1,\, a}
{\bar t}_{m b},\label{4.7}\\
d L_{E_{m+1, m}}(C_{a b}^{(r+1)})&=& -(-1)^{[a]+[b]} t_{m a}
{\bar t}_{m+1,  b},  \quad a, b \in{\bf I}. \eea

Note that as a $\U$-module, $\CG$ has a filtration defined by the
degrees of the polynomials in the $t_{ab}$ and the ${\bar
t}_{ab}$, and the filtration on the $\rm{U}(\fl_0)$-module
$\CG^{dL_{\rm{U}(\fl_0)}}$ defined by the degrees of the
polynomials in the $\{C_{ab}^{(i)} | i\ne r; a, b\in \bf{I}\}$ is
compatible with this filtration.  Thus in order to find those
$f\in \CG^{d L_{\rm{U}(\fl_0)}}$ such that $ d L_{E_{m m+1}}(f) =
L_{E_{m+1, m}}(f)=0$, by passing through to the associated graded
modules defined by these filtrations if necessary, we may assume
that $f$ is homogeneous of degree $\mu$ in the elements of
$\{C_{ab}^{(i)} | i\ne r; a, b\in \bf{I}\}$.  We consider an
element $f\in \CG^{d L_{\rm{U}(\fl_0)}}$ as a polynomial in
$\{C_{a b}^{(r+1)} | a, b\in \bf{I}\}$ with coefficients being
polynomials in $\{C_{a b}^{(i)} | i\ne r, r+1; \ a, b\in
\bf{I}\}$. Set $C_{ab} = C_{ab}^{(r+1)}$ ($a, b\in \bf{I}$).  Then
by Proposition 4.1, the subalgebra $<C_{ab} | a, b \in \bf{I}>$
has a basis consists of elements of the form \bea
C_{a_{1}b_{1}}^{p_{1}} \cdots C_{a_{s}b_{s}}^{p_{s}}
C_{c_{1}d_{1}} \cdots C_{c_{t}d_{t}}, \label{4.9}\eea with
$[a_{i}]+[b_{i}]=0$ ($1\le i\le s$), $[c_{j}]+[d_{j}]=1$ ($1\le
j\le t$), and $p_{i} \ge 0$ ($1\le i\le s$) are integers. Extend
such a basis of $<C_{ab} | a, b \in I>$ to a homogeneous basis
$\bf B$ of $\CG^{dL_{\rm{U}(\fl_0)}}$, so that the elements of
$\bf B$ are of the form \bea CC_{a_{1}b_{1}}^{p_{1}} \cdots
C_{a_{s}b_{s}}^{p_{s}} C_{c_{1}d_{1}} \cdots C_{c_{t}d_{t}},
\label{4.10}\eea where $C$ is a monomial in $\{C_{a b}^{(i)} |
i\ne r, r+1; \ a, b\in \bf{I}\}$. Now let us write $f = \sum_{0\le
k\le \mu} f_{k}$, where $f_{k}$ is a linear combination of the
basis elements of (\ref{4.10}) such that $\sum_{i}p_{i} + t = k$
and $\deg(C) + k = \deg(f)$. The action of $E_{m m+1}$ (similarly
for $E_{m+1 m}$) on the elements of (\ref{4.9}) can be computed by
using (\ref{4.7}), and we have \bea && d L_{E_{m m+1}}(
C_{a_{1}b_{1}}^{p_{1}} \cdots C_{a_{s}b_{s}}^{p_{s}}
C_{c_{1}d_{1}}
\cdots C_{c_{t}d_{t}}) = \nonumber\\
&&\quad -\sum_{i=1}^{s}(-1)^{t}p_{i} C_{a_{1}b_{1}}^{p_{1}}
\cdots C_{a_{i}b_{i}}^{p_{i}-1}
\cdots C_{a_{s}b_{s}}^{p_{s}} C_{c_{1}d_{1}}
\cdots C_{c_{t}d_{t}}t_{m+1 a_{i}}{\bar t}_{m b_{i}} \nonumber\\
&& \quad + C_{a_{1}b_{1}}^{p_{1}} \cdots C_{a_{s}b_{s}}^{p_{s}}
\sum_{j=1}^{t}(-1)^{j} C_{c_{1}d_{1}} \cdots {\hat
C}_{c_{j}d_{j}}\cdots C_{c_{t}d_{t}}t_{m+1 c_{j}}{\bar t}_{m
d_{j}}. \label{4.11}\eea Since the product map of $\CG$ is a
$\U$-module homomorphism (see \ref{modalg}), by (\ref{4.11}) the
action of $E_{m m+1}$ on the elements of (\ref{4.10}) is given by
\bea && d L_{E_{m m+1}}( CC_{a_{1}b_{1}}^{p_{1}} \cdots
C_{a_{s}b_{s}}^{p_{s}} C_{c_{1}d_{1}}
\cdots C_{c_{t}d_{t}}) = \nonumber\\
&&\quad -(-1)^{[C]}C\sum_{i=1}^{s}(-1)^{t}p_{i} C_{a_{1}b_{1}}^{p_{1}} \cdots
C_{a_{i}b_{i}}^{p_{i}-1}\cdots C_{a_{s}b_{s}}^{p_{s}} C_{c_{1}d_{1}} \cdots
C_{c_{t}d_{t}}t_{m+1 a_{i}}{\bar t}_{m b_{i}} \nonumber\\
&& \quad + (-1)^{[C]}CC_{a_{1}b_{1}}^{p_{1}} \cdots
C_{a_{s}b_{s}}^{p_{s}} \sum_{j=1}^{t}(-1)^{j}C_{c_{1}d_{1}} \cdots
{\hat C}_{c_{j}d_{j}}\cdots C_{c_{t}d_{t}}t_{m+1 c_{j}}{\bar t}_{m
d_{j}}, \label{4.12}\eea where ${\hat C}_{c_{j}d_{j}}$ means that
the factor $C_{c_{j} d_{j}}$ is omitted.

For an element $x$ of the form (\ref{4.10}), let $x'(ab)$ be \be C
C_{a_{1}b_{1}}^{p_{1}} \cdots C_{a_{i-1}b_{i-1}}^{p_{i-1}}
C_{a_{i}b_{i}}^{p_{i}-1} C_{a_{i+1}b_{i+1}}^{p_{i+1}} \cdots
C_{a_{s}b_{s}}^{p_{s}} C_{c_{1}d_{1}} \cdots C_{c_{t}d_{t}}, \ee
or \be CC_{a_{1}b_{1}}^{p_{1}} \cdots C_{a_{s}b_{s}}^{p_{s}}
C_{c_{1}d_{1}} \cdots {\hat C}_{c_{j}d_{j}}\cdots C_{c_{t}d_{t}},
\ee depending on whether $(ab)=(a_{i}b_{i})$ or
$(ab)=(c_{j}d_{j})$.

Let us make some observations. First note that since \be \langle
m_{\circ}(h\otimes  g),\, a\rangle = \langle h\otimes g, \,
\Delta(a)\rangle, \quad h,g\in\U^{0}, \quad a\in\U, \ee if $\{b_i |
1\le i\le \ell\}\subset \CG$ is a set of linearly independent
functions which are constants on $U(\fk_{0})$ and $\{t_{m+1 a}{\bar
t}_{m b} | a, b \in \bf{J}\subset \bf{I}\}$ is linearly independent,
then the set $\{b_{i}t_{m+1 a}{\bar t}_{m b} | 1\le i\le \ell, a, b
\in \bf{J} \}$ is linearly independent.  Then note that if $S_{cd}
\subset \bf{B}$ with $C_{cd}$ appearing in every element for a fixed
pair $c$ and $d$, then the set \be S'_{cd} = \{x'(cd)|x\in S_{cd}\} \ee
is linearly independent.  In fact the elements of $S_{cd}$ and
$C_{cd}S'_{cd}$ are the same up to signs.  Finally note that the only
relation amongst the elements in $\{t_{m+1 a}{\bar t}_{m b} | a, b \in
\bf{I}\}$ is (see (3.1)) \be \sum_{a\in I}t_{m+1 a}{\bar t}_{m
a}(-1)^{[a]} = \sum_{1\le a\le m} t_{m+1 a}{\bar t}_{m a} -
\sum_{m+1\le a\le m+n} t_{m+1 a}{\bar t}_{m a} = 0, \ee and this
relation can only come from (via the map $dL_{E_{m m+1}}$) \be
&& \sum_{1\le a\le m}C_{aa} - \sum_{m+1 \le a\le m+n} C_{aa}\\
&&\quad = \sum_{1\le a\le m} \sum_{c=m+1}^{l_{r+1}}t_{c a}{\bar t}_{c a}
- \sum_{m+1\le a\le m+n} \sum_{c=m+1}^{l_{r+1}}t_{c a}{\bar t}_{c a} = -k_{r+1},
\ee
i.e. a constant.

These observations together with (\ref{4.12}) imply that $ d
L_{E_{m m+1}}(f)= d L_{E_{m+1, m}}(f)=0$ if and only if $f =
f_{0}$, i.e. $f$ is independent of $C_{a b}^{(r+1)}$
($a,b\in\bf{I}$). Therefore,
\begin{lemma}\label{gCKG}
$\CKG$ is generated by the elements of \bea \{C_{a b}^{(i)} | i\ne r,
r+1; \ a, b\in {\bf I}\}.\label{CKGgen}\eea
\end{lemma}

By Theorem \ref{FFT} and the first fundamental theorem of the invariant
theory of the general linear group, $\CXbX^{\Psi(\Ul)}$ is generated
by ${\hat C}_{a b}^{(i)}$ ($i\ne r, r+1; a,b\in \bf{I}$), and ${\hat
C}_{a b}^{(r)}-{\hat C}_{a b}^{(r+1)}$ ($a,b\in \bf{I}$). We have
$\jmath({\hat C}_{a b}^{(i)})=C_{a b}^{(i)}$ ($i\ne r, r+1; a,b\in
\bf{I}$), which yield all the elements of (\ref{CKGgen}). This
establishes the short exact sequence \be &0\longrightarrow
\cJ^{\Psi(\Ul)} \longrightarrow \CXbX^{\Psi(\Ul)}{\longrightarrow}
\CKG\longrightarrow 0 \ee of $\Ul\otimes\Ul$-module algebras, thus
proves the first claim of Theorem \ref{construction}.

Let us now consider the subalgebra $\CKG^{d R_{\rm{U}(\fl_0)}}$ of
$\CKG$, which is generated by the elements of the set $\{ C^{(i, j)} |
i\ne r, r+1; j\ne r\}$, as follows from results of the last
subsection.  Amongst all the elements of this set, only $C^{(i, r+1)}$
are not annihilated by $d R_{E_{m m+1}}$ and $d R_{E_{m+1, m}}$. Thus
similar to the case of the left action, we can prove that $f\in \CKG^{d
R_{\rm{U}(\fl_0)}}$ satisfies $d R_{E_{m m+1}}(f)=0$ and $d R_{E_{m+1,
m}}(f)=0$ if and only if it is independent of the $C^{(i, r+1)}$
($i\ne r, r+1$). Observe that \be \CKGK =\{ f\in \CKG^{d
R_{\rm{U}(\fl_0)}} | d R_{E_{m m+1}}(f)=d R_{E_{m+1, m}}(f)=0\}. \ee
We have
\begin{lemma}\label{gCKGK}
$\CKGK $ is generated by the elements of \bea \{ C^{(i, j)} | i, j \ne
r, r+1\}.\label{CKGKgen}\eea
\end{lemma}

By Theorem \ref{FFT} and the first fundamental theorem of the
invariant theory of the general linear group,
$\CXbX^{\Psi(\Ul)\otimes\Phi(\Ul)}$ is generated by \be {\hat C}^{(i,
j)}, {\hat C}^{(r, j)}-{\hat C}^{(r+1, j)}, {\hat C}^{(i, r)}-{\hat
C}^{(i, r+1)}, \quad i, j\ne r, r+1, \quad a, b\in{\bf I}, \ee and
since \be \jmath(\{{\hat C}^{(i, j)}|i,j\ne r, r+1\})= \{{C}^{(i,
j)}|i, j\ne r, r+1\}, \ee we have the following short exact sequence
of $\Ul\otimes\Ul$-module algebras: \be 0\longrightarrow
\cJ^{\Psi(\Ul)\otimes\Phi(\Ul)} \longrightarrow
\CXbX^{\Psi(\Ul)\otimes\Phi(\Ul)}{\longrightarrow}
\CKGK\longrightarrow 0,\ee which is equivalent to the second claim of
Theorem \ref{construction}.

\begin{remark} \label{homogeneous-space} Geometric homogeneous
superspaces have been studied since the 1970s, see for example
\cite{Ko} and \cite{Ma88}.  Symmetric supersapces were also
classified by Serganova in \cite{Sergan} at the level of Lie
superalgebras. In relation to our algebraic definition of
homogeneous superspaces, one may ask the following question. Let
$P$ be the parabolic subgroup of $GL(m|n, \Lambda)$ with Lie
superalgebra $\fp$. We have the homogeneous superspace $GL(m|n,
\Lambda)/P$ (understood as a left coset of $P$). Now let $\fl$ be
the Levi factor of $\fp$ and take $\fk=\fl^{\sigma, \sqrt{i}}$,
with $\theta$ being the Hopf $\ast$-superalgebraic structure of
$\U$ corresponding to the compact real form of the general linear
superalgebra (defined by \eqref{compact-gl}). Then the question is
whether the homogeneous superspace determined by $\CKG$ is the
same as $GL(m|n, \Lambda)/P$ in some appropriate sense. We expect
the answer to be affirmative, but have not been able to locate a
reference, which addresses any form of the question, in the
literature on super-geometry.
\end{remark}

\section{Spherical Functions on $\CKG$ with Maximal Rank $K$}\label{max-rank}
We keep notations from the last section. In particular, we fix the
$\ast$-structure $\theta$ of $\U$ given by \eqref{compact-gl},
which corresponds to the real form $\fu(m|n)$ for the general
linear superalgebra. We use $\fl$ to denote the Levi factor of a
parabolic subalgebra of $\fg$, and set $\fk=\fl^{\sigma,
\sqrt{i}}$. The homogeneous superspaces studied in this section
are all examples of symmetric supersapces in the sense of
\cite{Sergan} (see Tables 2 and 3 in \cite{Sergan}).

\subsection{The case with $\fk=\fu(m|n-1)\oplus \fu(1) $}
We first examine in some detail the spherical functions on the
homogeneous superspace corresponding to $\fk=\fu(m|n-1)\oplus
\fu(1)$, where the complexification $\fl$ of $\fk$ is the
subalgebra of $\g$ spanned by the elements $E_{i j}$, $i, j\in
{\bf I}\backslash\{m+n\}$, and $E_{m+n\, m+n}$. But before
discussing the superalgebra $\CKG$, let us introduce the following
superalgebra.
\begin{definition} $\CS := \CG^{d L_{\rm{U}^\R(\fu(m|n-1)}}$
relative to $\fu(m|n-1)\subset\fk$. \end{definition} More explicitly,
\be \CS&=& \left\{f\in \CG \mid d L_k(f)=\epsilon(k) f, \ \forall k\in
\rm{U}^\R(\fu(m|n-1)) \right\}. \ee We can modify the analysis of
Subsection \ref{main} to construct $\CS$. With the help of Theorem
\ref{FFT} for ${\mathfrak{gl}}(m|n-1)$, we can show that $\CS$ is
generated by
\be z_a:=t_{m+n\, a},& & {\bar z}_a:={\bar t}_{m+n\, a},  \\
Q_{a b}&:=&\sum_{c<m+n}{\bar t}_{c a} t_{c b} (-1)^{[b][c]+[c]}, \quad a\in{\bf I},
\ee
where $z_a$ and ${\bar z}_a$ are odd if $a\le m$, and even otherwise.
The defining relations of $\CG$ imply $Q_{a b}=
\delta_{a b}\1- z_a {\bar z}_b (-1)^{[b]}.$ Thus
the $z_a$ and ${\bar z}_a$ generate $\CS$ by themselves. We have the
following result.
\begin{lemma} \label{StoP} The subalgebra of $\CS$ of $\CG$ is generated
by $z_a, \ {\bar z}_a$, $a\in{\bf I}$, which satisfy the following
relation \bea \sum_{a\in{\bf I}} {\bar z}_a z_a&=&\1. \label{sphere}
\eea
\end{lemma}

\begin{remark}
The notation suggests $\CS$  be the superalgebra of functions on the
supersphere. This can be understood as follows. Under the $\ast$-map
$\omega$ defined by \eqref{compact-G}, we have \be \omega(z_a)={\bar
z}_a, \quad \omega(\bar z_a)= z_a. \ee Thus we may interpret ${\bar
z}_a$ as the `complex conjugate' of $z_a$, and this indeed makes
perfect sense when $z_a$ and $\bar{z}_a$ are regarded as functions on
$G$ (see Subsection \ref{real}). Thus equation \eqref{sphere} defines
a supersphere in analogy with the embedding of a supersphere ${\mathbb
S}^{2n-1|2m}$ in $\C^{n|m}$. This also indicates the importance of the
$\ast$-structure in determining the underlying super manifold of
$\CKG$.
\end{remark}
\begin{remark}
When $\fk=\fu(m|n-1)\oplus\fu(1)$, we have $\CKG = \CS^{d
L_{\fu(1)}}$. This superalgebra embedding  $\CKG \hookrightarrow\CS$
corresponds to a projection from ${\mathbb S}^{2n-1|m}$ to the
symmetric superspace, which is the super generalization of the Hopf map
${\mathbb S}^{2n-1}\rightarrow {\mathbb C}{\mathbb P}^{n-1}$.
Therefore, we shall regard the symmetric superspace associated with
$\CKG$ as an algebraic analogue of the projective superspace ${\mathbb
C}{\mathbb P}^{n-1|m}$.
\end{remark}
We denote $\CKG$ by $\CP$ when $\fk=\fu(m|n-1)\oplus\fu(1)$. Lemma
\ref{StoP} immediately leads to the following result.
\begin{lemma}
The superalgebra $\CP$ is the $\ast$-subalgebra of $\CS$ generated by
the quadratic elements $z_a {\bar z}_b$, $a, b\in {\bf I}$.
\end{lemma}
\begin{proof}
Since for all $a$, $d L_{E_{m+n\, m+n}} z_a = z_a$, and $d L_{E_{m+n\,
m+n}} {\bar z}_a $ $=$  $ - {\bar z}_a$, any $d L_{\fu(1)}$-invariant
element of $\CS$ must be a polynomial in $z_a {\bar z}_b$, $a, b\in
{\bf I}$. This result can also be obtained in a more direct way by
using Theorem \ref{construction}.
\end{proof}
\begin{remark}
We should emphasize that elements of $\CP$ are functions on the
projective superspace that are  not `holomorphic' in general because
$\CP$ is a $\ast$-superalgebra.
\end{remark}

Now we use Theorem \ref{construction} to extract the algebra $\CP^{d
R_{\rm{U}^\R(\fk)}}$ of spherical functions on the projective
superspace. Let $z:=z_{m+n}$ and ${\bar z}={\bar z}_{m+n}.$ We have
\begin{theorem}\label{projective}
The algebra of the spherical functions on the projective superspace is
generated by $r:=z {\bar z}$ as a $\ast$-subalgebra of $\CP$. When
$n>1$, the spherical functions form a polynomial algebra in one
variable. When $n=1$, we have $(1-r)^{m+1}=0.$
\end{theorem}
\begin{proof}
It is an immediate consequence of Theorem \ref{construction} that the
algebra $\CP^{d R_{\rm{U}^\R(\fk)}}$ of the spherical functions on the
projective superspace is indeed generated by the single element $r$.

When $n=1$, all the $z_c$, ${\bar z}_c$, $c\le m$, are odd. Thus
the $(m+1)$-th power of $1-r=\sum_{c\le m} z_c {\bar z}_c$ vanishes
identically.

To study the case with $n>1$, we first analyse $\C[\rm{GL}_n]$, the
algebra generated by the matrix elements of the contravariant and
covariant vector representations of ${\mathfrak{gl}}(n)$. Let
$\fq={\mathfrak{gl}}(n-1)\oplus {\mathfrak{gl}}(1)$ be the subalgebra
of ${\mathfrak{gl}}(n)$ embedded block diagonally. Set $A=
\C[\rm{GL}_n]^{d L_{\rm{U}(\fq)}\otimes d R_{\rm{U}(\fq)}}$. Recall
that $\C[GL_n]$ is semi-simple as a left module $\rm{U}(\fq)$-module
under the action  ${d L_{\rm{U}(\fq)}\otimes d R_{\rm{U}(\fq)}}$. There
exists a surjective ${d L_{\rm{U}(\fq)}\otimes d
R_{\rm{U}(\fl)}}$-module map $\psi: \C[\rm{GL}_n] \rightarrow A$. Let
$\psi^*, (\id-\psi)^*: \rm{U}({\mathfrak{gl}}(n))\rightarrow
\rm{U}({\mathfrak{gl}}(n))$ be vector space maps defined by \be
\langle f,  (\id-\psi)^*(u)\rangle &=& \langle (\id-\psi)(f),
u\rangle,\\
\langle f,  \psi^*(u)\rangle &=& \langle \psi(f),  u\rangle, \quad
\forall u\in{\mathfrak{gl}}(n), f\in\C[\rm{GL}_n].\ee Since the dual
space pairing between $\C[\rm{GL}_n]$ and $\rm{U}({\mathfrak{gl}}(n))$
is non-degenerate, there is a non-degenerate pairing between $A$ and
$\psi^*(\rm{U}({\mathfrak{gl}}(n)))$.  Now as vector spaces, \be
\psi^*(\rm{U}({\mathfrak{gl}}(n)))&\cong&
\rm{U}({\mathfrak{gl}}(n))/\left(\fq
\rm{U}({\mathfrak{gl}}(n))+\rm{U}({\mathfrak{gl}}(n))\fq\right),\ee
where the right hand side is clearly infinite dimensional. This in
particular implies that the subalgebra $A$ of $\C[\rm{GL}_n]$ is
infinite dimensional.

Let $\zeta: \CP^{d R_{\rm{U}^\R(\fk)}}\rightarrow \C[\rm{GL}_n]$ be
the map defined for any $f\in \CP^{d R_{\rm{U}^\R(\fk)}}$ and
$u\in\rm{U}({\mathfrak{gl}}(n))$ by $\langle \zeta(f), u\rangle =
\langle \zeta(f), i(u)\rangle$, where $i$ is the canonical embedding
$\rm{U}({\mathfrak{gl}}(n))\subset \U$. Then $\zeta$ is an algebra
homomorphism, and we have \be \zeta(\CP^{d R_{\rm{U}^\R(\fk)}})&=&A.\ee

If there existed a non-trivial polynomial $P(r)$ in $r$ which was
identically zero as an element of $\CP$, then $\CP^{d
R_{\rm{U}^\R(\fk)}}$ would have to be finite dimensional over $\C$.
This contradicts the fact that $A$ is an infinite dimensional algebra.
\end{proof}

Let us now study the action of a generalized Laplacian operator on
the spherical functions. Recall that the quadratic Casimir of $\U$
can be expressed as $c=\sum_{a, b=1}^{m+n} (-1)^{[b]} E_{a b} E_{b
a}$. For any $f\in\CKGK$, we have $d R_X d R_c(f)$  $=$  $d R_c d
R_X(f)$  $=0,$ $\forall X\in\fl.$ That is $d R_c(f)\in\CKGK$.
Consider the following generalized Laplacian operator on the
homogeneous superspace: \be \nabla^2 &=&-\sum_{i=1}^{m+n-1} E_{i,
m+n} E_{m+n, i}. \ee Then the actions of $d R_{\nabla^2}$ and
$\frac{1}{2}d R_c$ coincide on $\CKGK$. Thus  $d R_{\nabla^2}$
also maps $\CKGK$ to itself.

In the case of the projective superspace, we can show that \bea d
R_{\nabla^2} (r^k) &=& k r^{k-1} \left[ (m-n-k+1) r + k\right], \quad
k=0, 1, .... \eea Let us now consider eigenfunctions of $d
R_{\nabla^2}$ in $\CP^{d R_{\rm{U}^\R(\fk)}}$. Things turn out to be
quite different for $m-n+1\le 0$ and $m-n+1>0$.
\begin{enumerate}
\item If $m-n+1\le 0$, there exists an eigenfunction $\theta_k\in
\CP^{d R_{\rm{U}^\R(\fk)}}$ of $d R_{\nabla^2}$ for each $k\in\Z_+$
with $d R_{\nabla^2}(\theta_k) = k (m-n-k+1) \theta_k, $ where \bea
\theta_k&=&\sum_{i=0}^k (-1)^i \left(\begin{array}{c} n-m+2k-2\\ i
\end{array} \right) \left(\begin{array}{c}
k\\ i \end{array} \right)^2 (i!)^2 r^{k-i}. \label{eigenfunction} \eea
Furthermore,  the $\theta_k$, $k\in\Z_+$, span $\CP^{d
R_{\rm{U}^\R(\fk)}}$.

\item If $m-n+1>0$, we let $L=m-n+1$, and denote by
$\left[\frac{L}{2}\right]$
the largest integer $\le L/2$. Then there exists an eigenfunction
$\theta_k\in\CP^{d R_{\rm{U}^\R(\fk)}}$ of $d R_{\nabla^2}$ for
each non-negative integer $k$ satisfying either $k\le
\left[\frac{L}{2}\right]$ or $k>L$ with $d
R_{\nabla^2}$-eigenvalue $k(L-k)$, where the $\theta_k$ are still
given by \eqref{eigenfunction}. However, the $\theta_k$'s do not
span $\CP^{d R_{\rm{U}^\R(\fk)}}$.
\end{enumerate}
Note that if $m-n+1>0$, the operator $d R_{\nabla^2}$ is not
daigonalizable over $\CKGK$.  The simplest illustration comes from
the case with $L=1$, where $\CKGK$ is the direct sum of $\{a+b
r|a, b\in \C\}$ and $\oplus_{k>1}\C\theta_k(r)$. While acting
diagonally on the latter subspace, $d R_{\nabla^2}$ acts on the
former subspace by $d R_{\nabla^2}(a+b r) = b$.

\begin{remark} $\CG$ is not semi-simple with respect to $d R_{\U}$.
There exist $d R_{\U}$-submodules of $\CG$ on which $d R_c$ can not be
diagonalized. Therefore, $d R_{\nabla^2}$ is not diagonalizable on
$\CKGK$ in general, and case (2) shows this fact.
\end{remark}

\subsection{The other maximal rank $K$ cases} We assume that both $m$ and $n$
are greater than $2$ in this subsection, and consider the maximal
rank $K$'s that correspond to the subalgebras $\fk_{n,k}:= \fl_{n,
k}^{\sigma, \sqrt{i}}$ and $\fk_{m,k}:=\fl_{m,k}^{\sigma,
\sqrt{i}}$, where \be \fl_{n, k}&=& \mathfrak{gl}(m
| n-k) \oplus \mathfrak{gl}(k), \quad 0<k\le n,\\
\fl_{m,k}&=& \mathfrak{gl}(k)\oplus \mathfrak{gl}(m-k | n), \quad
0<k\le m. \ee For the subalgebra $\fk_{n,k}$, by Theorem 4.1, the
corresponding homogenous superspace $\C(K_{n,k}\backslash G)$ is
generated by \be C_{ab} = \sum_{c=m+n-k+1}^{m+n}t_{ca}{\bar t}_{cb},
\quad a,b\in \bf{I}. \ee Note that $[c] = 1$. As in Theorem 5.1, we
can show that $\C[C_{ab}]$ forms a polynomial algebra in one variable
if $[a] = [b] =1$; and if $[a] = [b] = 0$, then $(C_{ab})^{k+1} = 0$
and $(C_{ab})^{k} \ne 0$. Recall that by Proposition 4.1, we always
have $(C_{ab})^{2} = 0$ if $[a] + [b] = 1$. The subalgebra of spherical
functions $\C[K_{n,k}\backslash G/K_{n,k}]$ is generated by \be C =
\sum_{c, a = m+n-k+1}^{m+n} t_{ca}{\bar t}_{ca}, \ee and forms a
polynomial algebra in one variable.  Similarly, for $\fk_{m,k}$, the
symmetric superspace $\C(K_{m,k}\backslash G)$ is generated by \be
C_{ab} = \sum_{c=1}^{k}t_{ca}{\bar t}_{cb}, \quad a,b\in \bf{I}. \ee
If $[a] = [b] = 0$, then $\C[C_{ab}]$ forms a polynomial algebra in one
variable, and if $[a] = [b] = 1$, then $(C_{ab})^{k+1} = 0$ and
$(C_{ab})^{k} \ne 0$. The subalgebra of spherical functions
$\C(K_{m,k}\backslash G/K_{m,k})$ is generated by \be C = \sum_{c, a
=1}^{k} t_{ca}{\bar t}_{ca}, \ee as a polynomial algebra.  To
summarize, we have
\begin{theorem} \label{maxrank}
1) If $m \le n$, then there is an onto algebra homomorphism \be \phi
:\C(K_{n,k}\backslash G) \rightarrow \C(K_{m,k}\backslash G) \ee which
induces an isomorphism $\C(K_{n,k}\backslash G/K_{n,k})\rightarrow
\C(K_{m,k}\backslash G/K_{m,k})$.

2) For each $1\le k < n$, there is an onto algebra homomorphism \be
\phi_{k+1,k} :\C(K_{n,k+1}\backslash G) \rightarrow
\C(K_{m,k}\backslash G) \ee which induces an isomorphism
$\C(K_{n,k+1}\backslash G/K_{n,k+1})\rightarrow \C(K_{n,k}\backslash
G/K_{n,k})$.
\end{theorem}
\begin{proof} For 1), we just need to note that any relation amongst
the $C_{ab}$ holds for both algebras by symmetry. For 2), let the
generators of $\C(K_{n,k}\backslash G)$ described above be $C_{ab}(k)$
($a,b\in {\bf I}, 1\le k \le n$), and define $\phi_{k+1,k}
:\C(K_{n,k+1}\backslash G) \rightarrow \C(K_{m,k}\backslash G)$ by
$\phi_{k+1,k}(C_{ab}(k+1)) = \frac{k+1}{k}C_{ab}(k)$.
\end{proof}

\section*{Acknowledgement}
We gratefully acknowledge financial support from the University of
Sydney and the Australian Research Council.

\end{document}